\documentclass{article}
\usepackage[utf8]{inputenc}

\usepackage{amssymb}
\usepackage{amsmath}
\usepackage{tikz}
\usepackage{tikz-cd}
\usepackage{enumitem}
\usepackage{accents}
\usetikzlibrary{cd}
\usetikzlibrary{babel}

\usepackage{biblatex}
\newcommand{\K}[1]{\mathcal{K}(\mathcal{#1})}

\newcommand{\D}[1]{\mathcal{D}(\mathcal{#1})}

\newcommand{\ChainComplex}[1]{{#1}_{\bullet}}
\newcommand{\CochainComplex}[1]{{#1}^{\bullet}}

\title{The Bornological Verdier Dual of the Structure Sheaves of Complex Manifolds}
\author{Christopher Burns}
\date{December 2022}

\begin{document}

\maketitle

\begin{center}
    \textbf{Abstract}
\end{center}

Quasi-abelian categories and their homological properties are recalled from \cite{schneiders1999quasi} in the first section, and the following section recalls the details of bornological sheaves from \cite{meyer2007local} and \cite{prosmans2000topological}. The third section provides an alternative proof of the bornological Verdier duality of \cite{prosmans2000topological}, using the theory of residues and duality utilised in \cite{chiarellotto1990duality} and \cite{van1992serre} in the rigid analytic context. 

\section{Quasi-abelian Categories and Sheaves}

\large

\noindent \textbf{Quasi-abelian categories}\\

\normalsize

Let $\mathcal{C}$ be an additive category with kernels and cokernels.\\

\noindent \textbf{Definition} - A morphism $f:X \rightarrow Y$ is \textit{strict} if it induces an isomorphism $\Tilde{f}: \mathrm{coim}(f) \rightarrow \mathrm{im}(f)$.\\

The category $\mathcal{C}$ is abelian if all maps are strict. Remark also that kernels and cokernels are always strict, and so any strict monomorphism admits a decomposition as a strict epimorphism, an isomorphism (which can simply be ignored by absorbing into on of the factors), and a strict monomorphism. Quasi-abelian categories instead satisfy a stability condition on strict morphisms, and so they generalise abelian categories:\\

\noindent \textbf{Definition} - An additive category $\mathcal{C}$ with kernels and cokernels is quasi abelian, if:

\begin{enumerate}

    \item The pullback of a strict epimorphism is strict: in the following pullback diagram, $f'$  strict implies that $f$ strict:
    
    \large
    
    \begin{center}
        \begin{tikzcd}
        X \arrow[r, "f"] \dar & Y \dar \\
        X' \arrow[r, "f'"] & Y'
        \end{tikzcd}
    \end{center}
    
    \normalsize
    
    \item Dually, the pushforward of a strict monomorphism is a strict monomorphism - so in the diagram above, instead take it to be a pushout square. Then $f$ a strict monomorphism implies $f'$ a strict monomorphism.

\end{enumerate}

\noindent \textbf{Examples} - Functional Analysis is a rich source of quasi-abelian categories, such as topological abelian groups, locally convex topological vector spaces and, in particular, Fréchet spaces. Of particular interest to us are complete convex bornological spaces - vector spaces with a distinguished class of \textit{bounded} subsets. This category will be denoted $\mathrm{CBorn}$.\\

Quasi-abelian categories possess familiar homological properties, which in turn also grants sheaves valued in them such familiar properties, while also being much more rigid in structure than, say, Quillen exact categories. We summarise those properties proven in \cite{schneiders1999quasi} - let $\mathcal{C}$ be a quasi-abelian category:

\begin{enumerate}

    \item We may place $t$-structures on the homotopy category $\K{E}$, but because the image and coimage of a map are generally distinct, there is no distinguished candidate to describe left truncation. This leads to two $t$-structures, the \textit{left and right} $t$-structures. All concepts dependent on the choice of $t$-structure, such as the hearts, will be called left or right accordingly.\\

    \item The \textit{derived category} $\D{E}$ is the localisation of $\K{E}$ with respect to the saturated null system of \textit{strictly exact complexes}. A complex is strictly exact in degree $n$, if $\mathrm{ker}(d^n) \cong \mathrm{im}(d^{n-1})$, and $d^n$ is strict, and a complex is strict if it's strict in all degrees. 
    
    \item Because homotopic complexes are strictly exact in the same degrees, and strict exactness in a given degree is inherited by mapping cones, it follows that the $t$-structures on $\K{E}$ localise to define $t$-structures on $\D{E}$.

    \item The left heart of $\D{E}$ is denoted $\mathcal{LH(E)}$, and is equivalent to the full subcategory of $\K{E}$ of complexes with a single nonzero monomorphic map between the only two nonzero terms. Taking the cokernels of such maps, and noting that the cokernels of parallel arrows of a cocartesian diagram are isomorphic, this enables us to embed $C$ in $\mathcal{LH(C)}$, and this embedding has a left adjoint. This embedding induces a derived isomorphism

    \begin{equation}
        \D{E} \cong \D{LH(E)}
    \end{equation}

    \noindent Therefore, at the derived level, quasi-abelian categories can be understood in terms of the conceptually easier abelian heart (recall that the heart of a t-structure is always an abelian category - see \cite{gelfand2002methods}).

    \item Derived functors can be defined similarly to the classical case as in \cite{gelfand2002methods}, with only minor technical adjustments.

    \item The category of sheaves on a topological space $X$ valued in quasi-abelian $C$, $\mathrm{Sh}(X; \mathcal{C})$, is quasi-abelian. Strict exactness can be measured stalkwise, and the category admits strong properties with respect to limits, generators and exactness, so that the operations used later in the category of complex convex bornological sheaves will be valid.

    \item $\mathcal{LH}(\mathrm{Sh}(X; \mathcal{C})) \cong \mathrm{Sh}(X; \mathcal{LH(C)})$, and combining this with the previous derived equivalence, we see that sheaves valued in quasi-abelian categories may also be understood in purely abelian terms at the derived level:

    \begin{equation}
        \mathcal{D}(\mathrm{Sh}(X; \mathcal{LH(C)})) \simeq \mathcal{D}(\mathrm{Sh}(X; \mathcal{C}))
    \end{equation}

    From this one can derive theorems for derived categories of sheaves valued in a quasi abelian category from analogous results for abelian categories.
    
\end{enumerate}

\newpage

\section{Ind(Ban) and Bornological Sheaves}

We are interested in convex bornological sheaves for two reasons. Firstly, most of the sheaves of interest in algebraic analysis are valued in some functional-analytic category which embeds into $\mathrm{CBorn}$. Secondly, $\mathrm{CBorn}$ is a particularly well behaved quasi-abelian category, with favourable closed structure. In \cite{prosmans2000topological}, results are proven for the category $\mathrm{Ind(Ban)}$, but in fact the subcategory of reduced inductive systems is equivalent to $\mathrm{CBorn}$ \cite{meyer2007local}.\\

\noindent \textbf{Definition} - The functor $\mathrm{IB}:\mathrm{Tc} \rightarrow \mathrm{Ind(Ban)}$ is defined on objects by

\begin{equation}
    \mathrm{IB}(V) = \underset{B \in \mathcal{B}_V}{``\varinjlim"} \hat{V}_B
\end{equation}

\normalsize

\noindent where $\mathcal{B}_V$ is the category of absolutely convex bounded subsets $B \subseteq V$, and $\hat{V}_B$ is the completion of the linear span of $B$. A continuous map of locally convex topological vector spaces will induce maps on all these linear spans, and hence on their completions. These maps extend to the formal inductive limits, thereby defining a functor $\mathrm{IB}:\mathrm{Tc} \rightarrow \mathrm{Ind(Ban)}.$ This may then be applied to the structure sheaf of a complex manifold, viewed as a sheaf valued in $\mathrm{Tc}$.\\

Let us now discuss how we can adapt this functor into one which maps to the conceptually simpler category of complete bornological vector spaces:\\

\noindent \textbf{Definitions} - A \textit{disc}, or a \textit{balanced or circled set} in a real or complex vector space, is convex, absolutely convex, and internally closed. Intuitively, a disc is a subset, which remains invariant under multiplication by the unit disc in $\mathbb{R}$ or $\mathbb{C}$.\\

To any disc, we may associate the gauge seminorm and the associated seminormed space $V_B$ on its linear span. The disc is \textit{complete} if $V_B$ is complete.\\

A bornological space is then complete if every bounded set is contained in a complete bounded disc.\\

Define an order on the discs of $V$: $B \leq B'$ if $B'$ \textit{absorbs} $B$ - that is, $B \subseteq rB'$ for some positive real number $r$. If $B \leq B'$, there is an injective bounded linear map $V_B \rightarrow V_{B'}$. Given a bounded linear map $V \rightarrow W$, by definition this induces a map of the directed sets of bounded discs, which defines a morphism of inductive systems. This functor is called $\textit{dissection}$ and denoted $\mathrm{diss}$.\\

On an ordinary bornological space, the target of $\mathrm{diss}$ is $\mathrm{Ind}(\mathrm{Sns})$, the Ind completion of the category of seminormed spaces. If $V$ is complete then, by assumption, the spaces $V_B$ are complete, and so $\mathrm{diss}$ may be regarded as a functor $\mathrm{CBorn} \rightarrow \mathrm{Ind}(\mathrm{Ban})$. The main theorem relevance to us is the following:\\

\noindent \textbf{Theorem} \cite{meyer2007local} - $\mathrm{diss}:\mathrm{CBorn} \rightarrow \mathrm{Ind}(\mathrm{Ban})$ is fully faithful, with essential image the \textit{reduced} systems - those with injective transition functions. We therefore have the fundamental equivalence:

\begin{equation}
    \mathrm{CBorn} \simeq \mathrm{Ind}(\mathrm{Ban})_{\mathrm{red}}
\end{equation}

In \cite{prosmans2000topological}, the functor IB associates to a locally convex topological vector space, a \textit{reduced} inductive system of Banach spaces. We may therefore restrict the codomain of IB and regard it as a functor valued in complete bornological vector spaces. By abuse of notation, we shall denote the associated functor by $\mathrm{IB}:\mathrm{Tc} \rightarrow \mathrm{CBorn}$.\\

We summarise the main facts about IB proven in \cite{prosmans2000topological}:

\begin{enumerate}
   \item $\mathrm{IB}$ commutes with \textit{reduced} inductive limits of systems of Fréchet spaces over $\mathbb{N}$ - those with all transition maps injective:

\begin{equation}
    \underset{n \in \mathbb{N}}{\varinjlim} \mathrm{IB}(F_n) \cong \mathrm{IB}(\underset{n \in \mathbb{N}}{\varinjlim} (F_n))
\end{equation}

\item For $V, W \in \mathrm{Tc}$ we can bestow the space $L_b(V, W) = \mathrm{Hom}_\mathrm{Tc}(V, W)$ with the structure of a locally convex topological vector space, with seminorms 

$$\{p_B: B \subseteq V \mathrm{bounded}, p \, \mathrm{a \, continuous \, seminorm \, on \, W}\}$$

\noindent and $p_B(f) = \underset{v \in B}{\mathrm{sup}}p(f(e))$. So, $p_B$ evaluates on $f$ the supremal value that $p \circ f$ takes on $B$.\\

The tensor product of two locally convex topological vector spaces $V, W$, inherits a canonical family of seminorms $\{p \otimes q\}$ for $p, q$ seminorms on $V, W$ respectively. For $x$ an element of $V \otimes W$

\begin{equation}
    (p \otimes q)(x) = \underset{x = \underset{i \in I}{\Sigma} v_i \otimes w_i}{\mathrm{inf}} \Sigma p(v_i)q(w_j)
\end{equation}

where the infimum runs over all possible representations of $x$ as an element of the tensor product. We record the isomorphisms of note - firstly for $V$ bornological and $W$ complete:

\begin{equation} \label{IB of bounded maps}
    \mathrm{IB}(L_b(V, W)) \cong \mathrm{Hom}_\mathrm{Ind(Ban)}(\mathrm{IB}(V), \mathrm{IB}(W))
\end{equation}

For arbitrary locally convex $V, W$:

\begin{equation}
    \mathrm{IB}(V) \hat{\otimes} \mathrm{IB}(W) \cong \mathrm{IB}(E \otimes F)
\end{equation}

where $E \otimes F$ has the locally convex structure from the system of seminorms defined above, and the internal tensor product on $\mathrm{Ind(Ban)}$ objects is obtained by the unique extension of the closed structure to Ind objects.

\item Whenever $E$ is a DFN space and $F$ is an FN space:

\begin{align}
    \mathrm{Hom}(\mathrm{IB}(E), \mathrm{IB}(F)) &\cong \mathbb{R}\mathrm{Hom}(\mathrm{IB}(E), \mathrm{IB}(F)) \\
    L(\mathrm{IB}(E), \mathrm{IB}(F)) &\cong \mathbb{R}L(\mathrm{IB}(E), \mathrm{IB}(F))
\end{align}

and whenever $X, Y$ objects of $\mathrm{Ind(Ban)}$ with $X$ nuclear (meaning that for any index $i$, there is a $j \geq i$ with $X(i) \rightarrow x(j)$ nuclear), then 

\begin{equation}
    E \hat{\otimes}^\mathbb{L} F \cong E \hat{\otimes} F
\end{equation}

\item IB respects sheaves in the following sense: if $F$ is a presheaf of Fréchet spaces, which is a sheaf when viewed as a presheaf of vector spaces on second countable space $X$ (such as a complex manifold), then the presheaf characterised by

\begin{equation} \label{IB Preserves Sheaves}
    \mathrm{IB}(F)(U) = \mathrm{IB}(F(U))
\end{equation}

is a sheaf.

\item Cartan's Theorem B holds for $\mathrm{IB}(\mathcal{O}_X)$, which by (\ref{IB Preserves Sheaves}) is a sheaf:

\begin{equation} \label{Cartan's Theorem B}
    \mathbb{R}\Gamma(U, \mathrm{IB}(\mathcal{O}_X)) \cong \mathrm{IB}(\mathcal{O}_X(U))
\end{equation}

\noindent whenever $U$ has vanishing algebraic nonzero cohomology cohomology, $H^k(U, \mathcal{O}_X) = 0, k > 0$.

\end{enumerate}

\newpage

\section{Verdier Dual of $\mathrm{IB}(\mathcal{O}_X)$}

In this chapter we present an alternative calculation of the Verdier dual of $\mathrm{IB}(\mathcal{O}_X)$ and $\mathrm{IB}(\Omega_X^{n-p})$, inspired by the argument from \cite{van1992serre}. By the equivalence between reduced inductive systems of Banach spaces and complete bornological vector spaces, we shall henceforth consider complete bornological sheaves on an $n$-dimensional complex manifold $X$. All derived functors are in the sense of Schneiders \cite{schneiders1999quasi}.\\

\noindent \textbf{Definition} - The \textit{Verdier dual} of a bornological sheaf $F$ on $X$, $D(F)$ is characterised by the formula, for all open $U \subseteq X$:\\

\begin{equation}
    D(F)(U) = \mathbb{R}L(\mathbb{R}\Gamma_c(U; F), \mathrm{IB}(\mathbb{C}))
\end{equation}

\noindent \textbf{Theorem} (Verdier Dual of $\mathrm{IB}(\mathcal{O}_X)$) - There are isomorphisms for all positive integers $p \leq n$:

\begin{equation} \label{Verdier Duality for }
    \mathrm{IB}(\Omega_X^{n-p}) \simeq D(\mathrm{IB}(\Omega_X^p))
\end{equation}

In \cite{prosmans2000topological} this result is proven by constructing a perfect pairing between $\mathrm{IB}(\Omega_X^{n-p})$ and $\mathrm{IB}(\Omega_X^p)$ by generalising the cup product and fibre integration to sheaves valued in $\mathrm{Ind(Ban)}$. Their argument thereafter is a double induction on the dimensions of the spaces considered and the number of irreducible components involved. After a long series of reductions, it is eventually shown to be a corollary of Stokes' theorem. \\

We propose a different proof, using an approach by descent to the calculation of the Verdier dual in the non-Archimedean context, detailed in \cite{van1992serre}. Rather than making the long series of reductions as in the above proof, the idea of this proof is as follows:

\begin{enumerate}
    \item Carry out the proof in the special case of complex affine space, $\mathbb{C}^n$. This is done by the construction of a \textit{residue map} on compactly supported cohomology, $H_c^n(X, \omega) \rightarrow \mathbb{C}$, which composed with a natural pairing on $H^0(\mathcal{O}_X) \times H_c^n(X, \omega)$ establishes duality.
    
    \item Given a closed immersion $i:X \rightarrow \mathbb{C}^n$ of a complex submanifold, we can express cohomology on $X$ in terms of cohomology of $\mathbb{C}^n$, using the argument given in \cite{chiarellotto1990duality}. Moreover, a residue map on $X$ can be defined in terms of the residue on the ambient affine space, independent of the choice of Stein neighbourhood and of the choice closed immersion $i:X \rightarrow \mathbb{C}^n$.
    
    \item Any \textit{Stein manifold} admits a closed immersion into some complex affine space. For details see \cite{grauert2013theory} or \cite{kaup2011holomorphic}. Therefore the duality theorem can be deduced for Stein domains from the affine case.
    
    \item Any point on a complex manifold admits a Stein neighbourhood, which may be assumed to be irreducible (and of course of constant dimension), and such that the intersections of these spaces are also irreducible Stein spaces. This is because holomorphic separability and convexity are preserved by the intersection of such domains, so one only has to take irreducible components.
    
    \item We therefore have unique canonical local duality isomorphisms, induced by unique canonical residue maps on these Stein neighbourhoods, and also on the overlaps of these domains. We conclude that all the duality isomorphisms glue to a global duality isomorphism for any complex manifold.
    
\end{enumerate}

(4) and (5) require no further explanation, so below we give a proof of (1), and briefly explain the pulling back of cohomology in (2) for closed immersed Stein domains. Before delving into the proof, we recall the necessary facts about Stein domains:\\

\noindent \textbf{Definition} \cite{kaup2011holomorphic} - A \textit{Stein space} $(X, \mathcal{O}_X)$ is a complex space, whose connected components are second countable as subspaces, which satisfies the following:

\begin{enumerate}
    \item $X$ is \textit{holomorphically separable} - for any distinct points $x,y$ of $X$, there is a global section $f \in \mathcal{O}_X$ with $f(x) \neq f(y)$.
    \item $X$ is \textit{holomorphically convex} - for any compact subset $K \subseteq X$, its \textit{holomorphically convex hull}
    
    \begin{equation}
        \hat{K}_{\mathcal{O}_X} = \{x \in X: \, \forall f \in \mathcal{O}_X, \, |f(x)| \leq \underset{k \in K}{\mathrm{sup}} |f(k)|\}
    \end{equation}
    
    \noindent is convex.
\end{enumerate}

Stein spaces are those complex spaces which possess very nice function theory. The criteria above guarantees the existence of sufficiently many global sections. Indeed, holomorphically convex subspaces of $\mathbb{C}^n$ are precisely the \textit{regions of holomorphy}: those spaces one which global holomorphic functions exist and cannot be analytically continued at any point. Holomorphic separability guarantees that there is no redundancy in the function theory.\\

Perhaps the most important theorem pertaining to Stein spaces is Cartan's Theorem B:\\

\noindent \textbf{Theorem} \cite{kaup2011holomorphic} - The analytic cohomology of a complex space $X$ is trivial if and only if $X$ is a Stein space.\\

Stein domains are therefore well suited for cohomological considerations. The following theorems demonstrate the flexibility of Stein manifolds for the local study of complex manifolds:\\

\noindent \textbf{Theorem} \cite{remmert1956espaces} - A Stein manifold $X$ of dimension $n$ admits a proper embedding into $\mathbb{C}^{2n+1}$.\\

\noindent \textbf{Theorem} - \cite{kaup2011holomorphic} - The topology of every complex space has a basis of open Stein neighbourhoods. In particular, any point of a complex manifold possesses a neighbourhood basis of Stein manifolds.\\

Therefore to locally study a complex manifold, we may take a Stein neighbourhood, which has nice cohomological properties, and exploit the cohomological relationship between the Stein manifold and a complex space into which it embeds. This is the approach we take in our alternative calculation of the Verdier dual. We need only a few additional facts:\\

\noindent \textbf{Theorem} (Extension Principle) \cite{grauert2012coherent} - For a closed complex subspace $X \subseteq Y$, an analytic sheaf $F$ on $X$ is coherent if and only if its pushforward $i_*F$ is coherent. This will apply in particular to the embeddings of Stein manifolds into affine space.\\

\noindent \textbf{Theorem} (Conormal Exact Sequence) \cite{hartshorne2013algebraic} - For an embedding of smooth complex manifolds $Y \hookrightarrow X$, if $\mathcal{I}$ is the associated ideal sheaf, there is a canonical exact sequence

\begin{center}
\begin{equation} \label{Conormal Sequence}
        \begin{tikzcd} 
    0 \rar & \mathcal{I}/{\mathcal{I}^2} \rar & {\Omega_X|_Y}^1 \rar & \Omega_Y^1 \rar & 0
\end{tikzcd}
\end{equation}

\end{center}

Moreover, this is an exact sequence of Fréchet spaces, and therefore an exact sequence of bornological sheaves for the canonically induced bornologies.\\

Through the use of Koszul resolutions, we can understand the cohomology of an embedded Stein manifold, or indeed any closed complex subspace, in terms of its local equations. First we need a definition:\\

\noindent \textbf{Definition}  \cite{kaup2011holomorphic} \cite{weibel1994introduction} - Let $R$ be a ring. A finite sequence of elements $(f_1, \cdots, f_n)$ of $R$ is a \textit{regular sequence} if each $f_i$ is a nonzerodivisor modulo $f_1, \cdots, f_{i-1}$.\\ 

We apply this definition to the local ring of a point (without loss of generality the origin) in complex affine space in the image of an embedded smooth Stein manifold $X$. We remark that the local equations $f_i$ for $X$ can be chosen such that the associated stalk elements $(f_i)_0$ form a regular sequence by smoothness - the Stein arises locally as an embedded subspace, and thus as the zero locus of complex variables. This could in principle fail for other embeddings, such as the embedding of a union of two spaces which meet tangentially at a point.\\

\noindent \textbf{Theorem} (Koszul Resolution)\cite{kaup2011holomorphic} - For any open region $U \subseteq \mathbb{C}^n$, and $i:X \hookrightarrow U$, $i(x) = 0$ the embedding of an analytic subset defined by equations $f_1, \cdots, f_p$, such that $(f_1)_0, \cdots, (f_p)_0$ defines a regular $(\mathcal{O}_U)_0$-sequence, the pushforward $i_*\mathcal{O}_X$ admits a free resolution locally at $x$, of length at most $p$. That is, on some open neighbourhood $0 \in V \subseteq U$, there is an exact sequence:

\begin{center}
    \begin{equation}
        \begin{tikzcd}
            0 \rar & \mathcal{O}_V^{m_p} \rar & \mathcal{O}_V^{m_{p-1}} \rar & \cdots \rar & \mathcal{O}_V^{m_0} \rar & (i_*\mathcal{O}_X)|_V \rar &0
        \end{tikzcd}
    \end{equation}
\end{center}

Of course, by translation, the assertion that $i(x) = 0$ is redundant by translation, so we may conclude that for any embedded Stein manifold into affine space, the structure sheaf admits a free resolution locally.\\

\newpage

\noindent \textbf{Proof of Verdier Duality in the Affine Case, $U = \mathbb{C}^n$}\\

\begin{enumerate}
    \item As in the proof in \cite{prosmans2000topological}, all the cases for various $p$ are equivalent, so we shall consider only the case of top differential forms.\\
    
    \item $\mathbb{C}^n$ possesses a filter by compact polydiscs $P_m$ of polyradius $(m, \cdots, m)$, for $m$ running over all positive integers. These polydiscs are cofinal in the system of compact subsets under inclusions, and so a compactly supported section is supported on one of these balls.
    
    \item By excision and Cartan's Theorem B for coherent sheaves, we have an isomorphism
    
    \begin{equation}
        H_{P_m}^{n}(U, \Omega_U) \cong H^{n-1}(U \backslash P_m, \Omega_U)
    \end{equation}
    
     \item $U \backslash P_m$ comes equipped with an acyclic cover by those subspaces, for varying $i$, with $i^{th}$ coordinate greater than $n$:
    
    \begin{equation}
        V_i = \{(z_1, ..., z_n): |z_i| \geq n \}
    \end{equation}
    
    These are Stein domains, and therefore defines an acyclic cover by Cartan's theorem B.  Therefore the associated Cech complex is adequate to calculate cohomology \cite{hartshorne2013algebraic}.
    
    \item Cech cohomology with respect to this cover, yields acyclicity in all degrees below $n$, the dimension of $U$. The $n^\mathrm{th}$ cohomology can be identified with the Fréchet space of Laurent series converging on the intersection of all of the acyclic domains
    
    \begin{center}
        \begin{equation}
        V = \{(z_1, ..., z_n): |z_i| \geq R_B \, \, \forall i \}
    \end{equation}
    \end{center}
    
    \noindent of the form
    
    \begin{equation} \label{Laurent Series Cohomology Space}
        \underset{m_i < 0}{\Sigma} a_{m}z^{m}dz/z
    \end{equation}
    
    \noindent where $m$ is the multi-index $(m_1, \cdots, m_d)$, $z^m = \prod_{i} z_i^{m_i}$, and $dz/z = \bigwedge_i {dz_i}/{z_i}$.\\
    
    \item As compactly supported cohomology is the colimit of cohomology supported on these balls, we find that $H_c^n(U, \Omega_U)$ can be identified with the space of such Laurent series, convergent outside \textit{some} ball $B$.
    
    \item Bearing in mind the statement of Verdier duality, we would like to relate this cohomology space to the space of bounded linear maps from $H^0(\mathcal{O}_X)$ to $\mathbb{C}$. Since all of the spaces involved are Fréchet spaces, boundedness and continuity are equivalent. Therefore, following the argument of \cite{van1992serre}, we associate such a continuous linear map to the Laurent series (\ref{Laurent Series Cohomology Space}) as follows: a global section in $H^0(\mathcal{O}_X)$ can be identified with a convergent power series, and we define the map by
    
    \begin{equation}
        f = \Sigma c_mz^n \mapsto \Sigma a_m c_m
    \end{equation}
    
    \noindent and since the projections onto the values of the coefficients are continuous and linear, every continuous linear map has this form. This proves the affine case:
    
    \begin{equation}
        H_c^n(\omega) \cong \mathrm{Ext}_c^n(\mathcal{O}_X, \omega)
    \end{equation}
    
    and the other Exts are zero, so by the above the dual is the translated structure sheaf.
    
    \item We may extend this result. Firstly, we have established a bilinear map, from which two isomorphisms can be deduced by dualising. Secondly we can consider more general scalars than those in $\mathbb{C}$, and consider scalars in an arbitrary $\mathbb{C}$-vector space $V$. We define the \textit{residue} $\mathrm{Res}: H_c^n(X, \omega) \rightarrow \mathbb{C}$ to be the map sending the class associated to a convergent sum (\ref{Laurent Series Cohomology Space}) to its lowest degree coefficient $a_0$. The \textit{trace} on $H_c^n(X, \omega \otimes V) \cong H_c^n(X, \omega) \otimes V$ is simply $\mathrm{Res} \otimes 1_V$. Then we deduce the duality results from the bilinear form
    
    \begin{align}
        H_c^n(X, \omega \otimes V) \times H^0(\mathcal{O}_X) &\rightarrow V \\
        \left(\underset{m_i < 0}{\Sigma} a_{m_i}z_i^{m_i}dz_1 \wedge dz_2 \wedge \cdots \wedge dz_d,  \Sigma c_mz^n \right) &\mapsto \Sigma a_m c_m
    \end{align}
    
    We deduce the isomorphisms
    
    \begin{align}
        H_c^n(\omega \otimes V) &\cong \mathrm{Ext}_c^n(\mathcal{O}_X, \omega \otimes V) \\
        H^0(\omega \otimes V) &\cong \mathrm{Ext}^0(\mathcal{O}_X, \omega \otimes V) 
    \end{align}
    
    \item Finally, these arguments can be generalised to a coherent $\mathcal{O}_X$-module $F$ in place of $\mathcal{O}_X$ by taking a local presentation. Therefore, we have the following isomorphisms for any coherent sheaf $F$ on $\mathbb{C}^n$:
    
    \begin{align}
        \mathrm{Ext}_c^p(F, \omega) &= 0, \, \, p \leq n\\
        \mathrm{Ext}_c^n(F, \omega) &\cong \mathrm{Hom \, cont}(H^0(F), \mathbb{C})\\
        \mathrm{Hom}(F, \omega) &\cong \mathrm{Hom \, cont}(H_c^n(F), \mathbb{C})
    \end{align}
    
    where $\mathrm{Hom \, cont}$ denotes the space of continuous, or equivalently bounded, linear maps between these spaces.\\
    
\end{enumerate}

\noindent \textbf{Remark} - Van der Put \cite{van1992serre} instead considers affine space over a rigid space $Y$, and so $V$ is replaced by a finitely generated $\mathcal{O}_Y$-module. For our purposes complex affine space is sufficient, and so $Y$ is taken to be a point, and $V$ is just a complex vector space.\\

\noindent \textbf{Pulling back Cohomology to an embedded smooth Stein manifold}\\

We follow the argument given in \cite{chiarellotto1990duality}. The argument is given in the rigid analytic context also, where Steins have a similar definition, and have similar properties to complex Stein domains:

\begin{enumerate}

    \item Let $i:X \hookrightarrow \mathbb{C}^n$ be a closed embedding, and cover the image of $X$ by open subsets $U_i$, such that on $U_i$, $X$ has local equations $f_1^i, \cdots, f_{n-q}^i$. These open sets may be chosen such that the stalks of the equations $f_j^i$ form a regular $\mathcal{O}_x$-sequence at any point $x \in X$. We can describe the structure sheaf on $X$ explicitly:
    
    \begin{equation} \label{Local defining equations for X}
    \mathcal{O}_X(V_i) \cong \mathcal{O}_{\mathbb{C}^n}(U_i)/(f_1^i, \cdots , f_{n-q}^i)
    \end{equation}

    \item Because $i$ is a closed immersion there is an \textit{underived} exceptional pullback $i^!$ - that is, a right adjoint to the pushforward $i_*$ on the category of sheaves. One finds by the following series of isomorphisms:
    
    \begin{align}
        \mathrm{Hom}_{\mathcal{O}_{\mathbb{C}^n}}(i_*\mathcal{F}, \mathcal{G}) &= \mathrm{Hom}_{i^{-1}\mathcal{O}_{\mathbb{C}^n}}(\mathcal{F}, i^!\mathcal{G}) \\
        &= \mathrm{Hom}_{i^{-1}\mathcal{O}_{\mathbb{C}^n}}(\mathcal{F}, \mathcal{H}\mathrm{om}_{i^{-1}\mathcal{O}_{\mathbb{C}^n}}(\mathcal{O}_X, i^!\mathcal{G}))\\
        &= \mathrm{Hom}_{i^{-1}\mathcal{O}_{\mathbb{C}^n}}(\mathcal{F}, i^*\mathcal{H}\mathrm{om}_{\mathcal{O}_{\mathbb{C}^n}}(i_*\mathcal{O}_X, \mathcal{G}))\\
        &= \mathrm{Hom}_{\mathcal{O}_X}(\mathcal{F}, i^*\mathcal{H}\mathrm{om}_{\mathcal{O}_{\mathbb{C}^n}}(i_*\mathcal{O}_X, \mathcal{G}))
    \end{align}
    
    Therefore by uniqueness of representatives we have a description of the functor $i^!$:
    
    \begin{equation} \label{Underived Exceptional Pullack}
        i^!\mathcal{G} \cong i^*\mathcal{H}\mathrm{om}_{\mathcal{O}_{\mathbb{C}^n}}(i_*\mathcal{O}_X, \mathcal{G})
    \end{equation}

    \item By taking the $n^\mathrm{th}$ exterior power, we get an isomorphism
   
    \begin{equation}
     i^*\Omega_{\mathbb{C}^n}^n \otimes \mathcal{O}_X \cong \bigwedge^{n-q} \mathcal{I}/{\mathcal{I}^2} \otimes \Omega_X^q
    \end{equation}
    
    This isomorphism is defined in the obvious way, as described in \cite{hartshorne1966residues}. That is, we can choose local generators $e_1, ..., e_n$ of $\Omega_{\mathbb{C}^n}^1 \otimes \mathcal{O}_X$, such that $e_1, ..., e_{n-q}$ generate $\mathcal{I}/{\mathcal{I}^2}$, and the images of the remaining generators generate $\Omega_X^1$. This choice determines the isomorphism.
    
    The wedge product on the right is invertible, yielding the isomorphism:
    
    \begin{equation} \label{Forms on X Equation}
        \Omega_X^q \cong i^*\Omega_{\mathbb{C}^n}^n \otimes \left(\bigwedge^{n-q} \mathcal{I}/{\mathcal{I}^2}\right)^*
    \end{equation}

    \item To analyse $\mathbb{R}i^!$ one considers the restriction to open affine neighbourhoods as in (\ref{Local defining equations for X}). Then, one can analyse the local Ext sheaves $\mathcal{E}\mathrm{xt}_{\mathcal{O}_X|_{U_i}}^j(\mathcal{O}_X|_{U_i}, \Omega^n_{\mathbb{C}^n}|_{U_i})$ using Koszul resolutions.
    
        \begin{enumerate}
        \item Consider a neighbourhood $U_i$ as in (\ref{Local defining equations for X}), so that $X$ is locally defined by the equations $f_1^i = \cdots = f_{n-q}^i = 0$.
        
        \item Denote the Koszul resolution of $i_*\mathcal{O}_X|_{U_i}$ by $\ChainComplex{K}(f_1^i, \cdots, f_{n-q}^i)$, and for any sheaf $F$ of $\mathcal{O}_{U_i}$-modules we define the \textit{dual complex}:
        
        \begin{equation}
            \CochainComplex{K}(f_1^i, \cdots, f_{n-q}^i; F) = \mathcal{H}\CochainComplex{\mathrm{om}}_{\mathcal{O}_X}(\ChainComplex{K}(f_1^i, \cdots, f_{n-q}^i), F)
        \end{equation}
        
        \noindent We see that because the Koszul complex is locally free, this calculates the sheaf Ext complex:
        
        \begin{align}
            \CochainComplex{K}(f_1^i, \cdots, f_{n-q}^i; F) &\cong \CochainComplex{\mathcal{E}\mathrm{xt}}_{\mathcal{O}_{U_i}}(\ChainComplex{K}(f_1^i, \cdots, f_{n-q}^i), F)) \\
            &\cong \CochainComplex{\mathcal{E}\mathrm{xt}}_{\mathcal{O}_{U_i}}(\mathcal{O}_{U_i}/(f_1^i, \cdots, f_{n-q}^i), F)) \\
            &\cong \CochainComplex{\mathcal{E}\mathrm{xt}}_{\mathcal{O}_{U_i}}(\mathcal{O}_X|_{V_i}, F))
        \end{align}
        
        We deduce an isomorphism on global cohomology
        
        \begin{equation} \label{n-q Ext isomorphism}
            \mathrm{Ext}^{n-q}(\mathcal{O}_X|_{V_i}, F) \rightarrow H^n\CochainComplex{K}(f_1^i, \cdots, f_{n-q}^i; F) \rightarrow  F/(f_1^i, \cdots, f_{n-q}^i)F
        \end{equation}
        
        The first map is the induced map on $(n-q)^\mathrm{th}$ cohomology. The second morphism $H^n(f_1^i, \cdots, f_{n-q}^i; F) \rightarrow F/(f_1^i, \cdots, f_{n-q}^i)F$, sends $\alpha \in K^n(f_1^i, \cdots, f_{n-q}^i; F)$ to $\alpha_{1, \cdots, n}$.
        
        \item The argument of Proposition 7.2 from \cite{hartshorne1966residues} adapts to our context, and allows us to conclude that
    
    \begin{equation}
        \mathcal{E}\mathrm{xt}_{\mathcal{O}_X|_{U_i}}^j(\mathcal{O}_X|_{U_i}, \Omega^n_{\mathbb{C}^n}|_{U_i}) = 0, \, \, j \neq n-q
    \end{equation}
    
    and therefore these Ext sheaves vanish globally also.

    \item By (\ref{n-q Ext isomorphism}) there are isomorphisms depending on the chosen equations for $X$
    
    \begin{equation}
        \phi_i: \mathcal{E}\mathrm{xt}_{\mathcal{O}_X|_{U_i}}^{n-q}(\mathcal{O}_X|_{U_i}, \Omega^n_{\mathbb{C}^n}|_{U_i}) \cong \frac{\Omega^n_{\mathbb{C}^n}|_{U_i}}{{\mathcal{I}|_{U_i}\Omega^n_{\mathbb{C}^n}|_{U_i}}}
    \end{equation}
    
    and these isomorphisms glue via the maps relating the systems of equations on overlaps. That is, given two systems of equations $X$ on $U_i$ and $U_j$, $(f_1^i, \cdots, f_{n-q}^i)$ and $(g_1^i, \cdots, g_{n-q}^i)$ respectively, then on the overlap $U_i \cap U_j$, these equations can be written in terms of each other. Equivalently, there is a matrix $c_{kl}$ with $g_k = \Sigma c_{kl}f_l$. By lemma 7.1 of \cite{hartshorne1966residues} these isomorphisms differ exactly by the canonical map $\mathrm{det}(c_{ij})$. These morphisms $\phi_i$ can therefore be glued along these maps, yielding a global isomorphism:
    
    \begin{equation}
        \phi:     \mathcal{E}\mathrm{xt}_{\mathcal{O}_X|_{U_i}}^{n-q}(\mathcal{O}_X|_{U_i}, \Omega^n_{\mathbb{C}^n}|_{U_i}) \cong \mathcal{H}\mathrm{om}_{\mathcal{O}_X}\left(\bigwedge^{n-q} \frac{\mathcal{I}}{\mathcal{I}^{2}} , \frac{\Omega^n_{\mathbb{C}^n}|_{U_i}}{{\mathcal{I}|_{U_i}\Omega^n_{\mathbb{C}^n}|_{U_i}}}\right)
    \end{equation}
    
    \noindent as the gluing data $\{\mathrm{det}(c_{ij})\}$ is described by the automorphisms of $\bigwedge^{n-q}\frac{\mathcal{I}}{\mathcal{I}^{2}}$. To be precise, this map $\phi$ evaluated on $f_1^i, \cdots, f_{n-q}^i$ is exactly $\phi_i$.
    
    \item Applying $i^*$ and the isomorphism (\ref{Forms on X Equation}), we deduce
    
    \begin{equation}
        i^*\mathcal{E}\mathrm{xt}_{\mathcal{O}_X}^i(\mathcal{O}_X, \Omega^n_{\mathbb{C}^n}) \cong i^*\Omega_{\mathbb{C}^n}^n \otimes \left(\bigwedge^{n-q} \mathcal{I}/{\mathcal{I}^2}\right)^* \cong \Omega_X^q
    \end{equation}
    
    We therefore finally arrive at a neat description of the derived exceptional pullback, according to (\ref{Underived Exceptional Pullack})
    
    \begin{equation}
        \mathbb{R}i^!(\Omega_{\mathbb{C}^n}^n) \cong \Omega_X^q[q-n]
    \end{equation}
    
    \end{enumerate}
    
    \item Recall that $i_*$ preserves coherence. To prove coherent duality we take an injective resolution of top forms
    
    \begin{center}
    \begin{tikzcd}
    0 \rar & \Omega_{\mathbb{C}^n}^n \rar & \CochainComplex{I}
    \end{tikzcd} 
    \end{center}

    \noindent and take cohomology of the adjunction isomorphism between $i_*$ and $i^!$. This yields the duality isomorphisms, for any coherent $\mathcal{O}_X$-module
    
    \begin{align}
        \mathrm{Ext}^i_{\mathcal{O}_{\mathbb{C}^n}}(i_*M, \Omega^n_{\mathbb{C}^n}) &\cong \mathrm{Ext}^{i-(n-q)}_{\mathcal{O}_X}(M, \Omega^q_{X}) \\
        \mathcal{E}\mathrm{xt}^i_{\mathcal{O}_{\mathbb{C}^n}}(i_*M, \Omega^n_{\mathbb{C}^n}) &\cong i_*\mathcal{E}\mathrm{xt}^{i-(n-q)}_{\mathcal{O}_X}(M, \Omega^q_{X})
    \end{align}
    
    \item Because $i$ is proper, the correspondence between compact subsets yields the isomorphism on compactly supported cohomology, for any sheaf of $\mathcal{O}_X$-modules $M$:
    
    \begin{equation}
        H_c^*(\mathbb{C}^n, i_*M) \cong H_c^*(X, M)
    \end{equation}
    
    These isomorphisms are already enough to pull back cohomology from affine space to the closed immersed Stein. From these we deduce the necessary isomorphisms. However, the isomorphisms need to be sufficiently canonical in order for us to deduce a global isomorphism. In particular, the isomorphism must not depend on the choice of embedding into affine space, and two different Stein covers should yield the same isomorphism. This is done by constructing an invariant trace map for Steins, based on the affine definition.
    
    \item We can define for an embedding $\phi:X \rightarrow \mathbb{C}^n$, a residue map on $X$:
    
    \begin{center}
        \begin{equation} \label{Residue on a Stein Manifold}
        \begin{tikzcd}[column sep = small]
             H_c^n(X, \Omega_X) \arrow[r, "\sim"] & H_c^n(\mathbb{C}^n, \mathcal{E}\mathrm{xt}^{n-q}(\phi_*\mathcal{O}_X, \Omega_{\mathbb{C}^n})) \arrow[r, "\alpha"] & \mathrm{Ext}_c^N(\phi_*\mathcal{O}_X, \Omega_{\mathbb{C}^n}) \arrow[dl, "\beta"] \\ & \mathrm{Hom \, cont}(H^0(\mathcal{O}_X), \mathbb{C}) \arrow[r, "\gamma"] & \mathbb{C} 
            \end{tikzcd}
         \end{equation}
    \end{center}
    
    The map $\alpha$ is an isomorphism, coming from the spectral sequence for $\mathrm{Ext}_c^N$. The isomorphism $\beta$ comes from the coherent duality isomorphism. The final map $\gamma$ is the evaluation on the constant function $1$ on $X$. The trace map is defined by taking the same sequence, tensoring the sheaves of differential forms (including the trivial one $\mathbb{C}$ on the point) with a complex vector space $V$.\\

    \item It remains to show the following, which we explain below, following the argument of \cite{van1992serre}: 
    
    \begin{enumerate}
        \item Given two distinct immersions into distinct affine spaces, the residues associated both equal the residue associated to the sum of the maps $(i_1, i_2):X \rightarrow \mathbb{C}^{(n_1 + n_2)}$ 
        
        \item Given an open immersion $X_1 \rightarrow X_2$ of Steins, the composition
    
    \begin{equation}
        H_c^n(X_1, \Omega_{X_1}) \rightarrow H_c^n(X_2, \Omega_{X_2}) \xrightarrow{\mathrm{res}_{X_2}} \mathbb{C}
    \end{equation}
    
    \noindent is $\mathrm{Res}_{X_1}$
    
    \end{enumerate}
\end{enumerate}

\large

\noindent \textbf{Invariance of Trace and Residue}\\

\normalsize 

\noindent \textbf{Lemma} - The residue map on affine space can be written as a composition

\begin{center}
    \begin{tikzcd}
    H_c^n(\mathbb{C}^n, \omega) \arrow[r, "\phi"] & H_c^n(\mathbb{P}^n, \omega) \arrow[r, "\sim"] & H^n(\mathbb{P}^n, \omega) \arrow[r, "\psi"] & \mathbb{C}
\end{tikzcd}
\end{center}

\noindent for some final isomorphism $\psi:H^n(\mathbb{P}^n, \omega) \rightarrow \mathbb{C}$.\\

\noindent \textbf{Proof} - by GAGA \cite{serre1956geometrie} the cohomology group of $\mathbb{P}^n$ can be calculated algebraically. To prove the statement we show that $\phi$ and $\mathrm{Res}_n$ have the same kernel. Then, $\psi$ can be chosen such that $\psi \circ \phi = \mathrm{Res}_n$.\\

We consider the action of $(S^1)^n$ on $\mathbb{C}^n$. It respects the polydiscs $P_m$, as well as the open cover of $\mathbb{C}^n \setminus P_m$, and therefore defines an action on $H_c^n(X, \omega)$. The extension of this action to $H^n(\mathbb{P}^n)$ is trivial. Therefore for any cohomology class $\xi$, and any operator $\sigma = (\lambda_1, \cdots, \lambda_n) \in (S^1)^n$ we have $\phi(\sigma\xi - \xi) = 0$. Representing $\xi$ as a convergent sum $\underset{m_i < 0}{\Sigma} a_{m}z^{m}dz/z$, we see that 

\begin{equation}
    \sigma\xi - \xi = \underset{m_i < 0}{\Sigma} a_{m}z^{m}(\lambda^m - 1)dz/z
\end{equation}

With careful choices of $\lambda$ such that the terms $\lambda^m-1$ don't vanish for $m \neq 0$, we see that the terms that vanish under $\phi$ are exactly those terms with $a_0 = 0$ - that is, $\mathrm{ker} \, \mathrm{Res}_n = \mathrm{ker} \phi$. This proves the factorisation, and then because $H^n(\mathbb{P}^n, \omega)$ is invariant under the projective general linear group, so is the residue.\\

\textbf{Lemma} - For any automorphism $\sigma$ of $\mathbb{C}^n$, there exists a unique invertible $f \in H^0(\mathcal{O}_X)$ such that $\mathrm{Res}_n \circ \sigma = f \cdot \mathrm{Res}_n$ as maps $H_c^n(\mathbb{C}^n, \omega) \rightarrow \mathbb{C}$.\\

\noindent \textbf{Proof} - The map $\mathrm{Res}_n \circ \sigma$ is a continuous linear map $H_c^n(\omega) \rightarrow \mathbb{C}$. Because $\omega$ is coherent, the duality statement implies that $\mathrm{Res}_n \circ \sigma(\xi) = \mathrm{Res}_n(f\xi)$. We may apply the same reasoning to $\sigma^{-1}$, and so we conclude that $f$ is invertible.\\

With this, we may conclude invariance of $\mathrm{Res}_n$ with respect to certain automorphisms of $\mathbb{C}^n$:\\

\noindent \textbf{Lemma} - Automorphisms of the form

\begin{equation}
    \sigma(z_1, \cdots, z_n) = (z_1, \cdots, z_m, z_{m+1} + k_{m+1}, \cdots, z_n + k_{n})
\end{equation}

\noindent where the $k_j$ are holomorphic functions varying with $z_1, \cdots, z_m$, are $\mathrm{Res}_n$-invariant: $\mathrm{Res}_n \circ \sigma = \mathrm{Res}_n$.\\

\noindent \textbf{Proof} - Such automorphisms can be written as a commutator composed with a translation. Commutators are $\mathrm{Res}_n$-invariant because of the last lemma. Translations are $\mathrm{Res}_n$-invariant from the power series formulation.\\

These lemmas are enough to prove the invariance properties of the generalised residue map (\ref{Residue on a Stein Manifold}).\\

\noindent \textbf{Theorem} - The residue (\ref{Residue on a Stein Manifold}) does not depend on the choice of embedding $\phi$.\\

\noindent \textbf{Proof} - A closed immersion of steins $\phi:X_1 \rightarrow X_2$ of dimensions $m_1, m_2$ respectively, induces a map $\Tilde{\phi}$ functorially on cohomology, by pulling back sheaves of differential forms:

\begin{center}
    \begin{tikzcd}[column sep=small]
    H_c^{m_1}(X_1, \Omega_{X_1}) \arrow[r, "\sim"] & H_c^{m_1}(X_2, \mathcal{E}\mathrm{xt}^{m_2-m_1}(\phi_*\mathcal{O}_X, \Omega_{X_2})) \rar & \mathrm{Ext}_c^{m_2}(\phi_*\mathcal{O}_{X_1}, \Omega_{X_2}) \dar \\ \, & \, & H_c^{m_2}(X_2, \Omega_{X_2})
    \end{tikzcd}
\end{center}

Now given two closed immersions $X \hookrightarrow \mathbb{C}^{n_1}$, $X \hookrightarrow \mathbb{C}^{n_2}$, we show that the residue associated to these embeddings is equal to the residue associated to the coproduct of the embeddings. Consider the commutative diagram

\large

\begin{center}

    \begin{tikzcd}
     & \mathbb{C}^{n_1} \arrow[rr, "l_1"] \arrow[dr, "k_1"] &  & \mathbb{C}^{n_1 + n_2} \arrow[dl, "\tau_1"]\\
    X \arrow[ur, "\phi_1"] \arrow[dr, "\phi_2", swap] \arrow[rr] &  & \mathbb{C}^{n_1 + n_2} & \\
     & \mathbb{C}^{n_2} \arrow[ur, "k_2"] \arrow[rr, "l_2"] &  & \mathbb{C}^{n_1 + n_2} \arrow[ul, "\tau_2"]
    \end{tikzcd}
    
\end{center}

\normalsize

\noindent where $k_1$ is the holomorphic extension of $\phi_2 \circ \phi_1^{-1}$ on $\phi_1(X)$, $l_1$ is the embedding $z \mapsto (z, 0)$, and this determines a unique map $\tau_1$. The other maps are defined similarly.\\

For the map $l_1$, the induced map $\Tilde{l_1}$ on cohomology induces the canonical map on power series rings. Therefore, it respects the lowest order coefficient, which is exactly the residue.\\

So, the residue associated with $\phi_1$ is the same as the residue associated with $l_1 \circ \phi_1$. As $\tau_1$ is an automorphism of the form considered before it respects residues. As this composition equals $k_1 \circ \phi_1$, the associated residues are equal. The same reasoning applies to the bottom half of the diagram, concluding the proof.\\

With this proven, we are entitled to unambiguously denote the residue on a Stein manifold $X$ by $\mathrm{Res}_X$, without reference to an ambient affine space into which it embeds. It remains only to prove the following:\\

\noindent \textbf{Proposition} - For a closed immersion of Stein spaces $\phi:X_1 \rightarrow X_2$, the composition

\begin{center}
    \begin{tikzcd}
    H_c^n(X_1, \Omega_{X_1}) \arrow[r, "\Tilde{\phi}"] & H_c^n(X_2, \Omega_{X_2}) \arrow[r, "\mathrm{Res}_{X_2}"] & \mathbb{C}
    \end{tikzcd}
\end{center}

\noindent is $\mathrm{Res}_{X_1}$. Therefore, the residue map and the coherent duality isomorphisms do not depend on the choice of Stein neighbourhood chosen.\\

\noindent \textbf{Proof} - $\mathrm{Res}_{X_2} \circ H_c^n(\phi)$ is continuous and linear, and so by coherent duality, it corresponds to a unique element $a(X_1, X_2) \in H^0(X_1, \mathcal{O}_{X_1})$:

\begin{equation}
    \mathrm{Res}_{X_2} \circ H_c^n(\phi)(\xi) = \mathrm{Res}_{X_1}(a(X_1, X_2)\xi)
\end{equation}

Our task is to show that $a(X_1, X_2) = 1$. First this is proven for the prototypical embedding of a polydisc into affine space:
$P^n(R) \subseteq \mathbb{C}^n$. The action of $(S^1)^n$ respects both the polydisc and affine space, and the residue maps are invariant under the action of each. From this symmetry, we see that the element $a(X_1, X_2)$ must be a scalar, because it's holomorphic and constant on tori.\\

Next we reduce to the one-dimensional case. We consider the commutative diagram:

\begin{center}
    \begin{tikzcd}
    P^n(R) \rar & \mathbb{C}^n \\
    P^1(R) \uar \rar & \mathbb{C} \uar
    \end{tikzcd}
\end{center}

\noindent where the vertical maps are the closed immersions, mapping to the first coordinate. Because these are closed immersions, the maps on cohomology are canonical, and the scalars $a(P^n(R), \mathbb{C}^n)$ and $a(P^1(R), \mathbb{C})$ must agree, reducing to dimension 1.\\

Van der Put very elegantly proves that $a^2 = a$ by considering the following diagram:

\begin{center}
    \begin{tikzcd}[column sep=small]
     & P^2(R) \arrow[rr] & & P^1(R) \times \mathbb{C} \arrow[rr] & & \mathbb{C}^2 &  \\
    0 \times P^1(R) \arrow[ur] \arrow[rr] & & 0 \times \mathbb{C} \arrow[ur] & & P^1(R)  \times 0 \arrow[ul] \arrow[rr] & & \mathbb{C} \times 0 \arrow[ul]
    \end{tikzcd}
\end{center}

Every map in this diagram is the evident embedding. Now on one hand, $a$ is determined by the composition of the top two horizontal maps, since this is the embedding of a polydisc. On the other hand, $a$ is also determined by both the horizontal lower maps, and since the upwards maps are all closed immersions, we conclude that $a^2 = a$, and $a \neq 0$ because $\mathrm{Res}_n \circ H_c^n(\phi) \neq 0$.\\

It remains to prove the result for any open immersion of Stein manifolds $\phi:X_1 \rightarrow X_2$ - that is, we'd like to show that for any point $p \in X_1$ that $a(X_1, X_2)(p) = 1$. Embed each $X_i$ into affine space, $\phi_i:X_i \rightarrow \mathbb{C}^{n_i}$, such that $p \mapsto 0$ under both embeddings. Choose small polydiscs $B_i \subseteq \mathbb{C}^{n_i}$ whose preimages in $X_1$ and $X_2$ are the same open neighbourhood $U$ of $p$ in $X_1$: $U = \phi_1^{-1}(B_1) = \phi_2^{-1}(B_2)$. The commutative diagrams

\begin{center}
    \begin{tikzcd}
    U \rar \dar & X_i \dar \\
    B_i \rar & \mathbb{C}^n_i
    \end{tikzcd}
\end{center}

\noindent demonstrate that $a(U, X_i) = 1$ as functions on $U$, and by functoriality of associated maps on cohomology we have
$a(U, X_2) = a(U, X_1)a(X_1, X_2)$, meaning that the restriction of $a(X_1, X_2)$ to any open neighbourhood is identically 1, and so it is globally constant of value 1, as required.\\

\newpage

\printbibliography

\end{document}